\documentclass[12pt,numbered,print,index]{amsart}
\usepackage[utf8]{inputenc}
\pdfoutput=1
\usepackage{adjustbox,amssymb,bbm,booktabs,enumitem,lmodern,mathrsfs,mathtools,microtype,tikz-cd,xcolor}
\usepackage{comment}
\usepackage[mathcal]{euscript}
\setenumerate{label={\normalfont(\arabic*)}}
\numberwithin{equation}{section}
\usepackage[T1]{fontenc}
\usepackage{color}
\usepackage{array}
\usepackage[backend=bibtex8,style=trad-alpha,doi=false,isbn=false,giveninits=true,mincitenames=3,maxcitenames=4]{biblatex}

\usepackage{setspace} 

\expandafter\let\csname ver@amsthm.sty\endcsname\relax
\let\theoremstyle\relax

\usepackage{amsthm}

\newtheorem{theo}{Theorem}[section]
\newtheorem*{theo*}{Theorem}
\newtheorem{cor}[theo]{Corollary} 
\newtheorem{lem}[theo]{Lemma}
\newtheorem{prop}[theo]{Proposition} 

\theoremstyle{definition}
\newtheorem{deff}[theo]{Definition}
\newtheorem{ex}[theo]{Example}
\newtheorem{ex*}{Example}
\newtheorem{rem}[theo]{Remark}

\DeclareMathOperator{\Lim}{Lim}
\DeclareMathOperator{\Liminf}{Lim\text{ }inf}
\DeclareMathOperator{\supp}{supp}
\DeclareMathOperator{\WAP}{\it WAP}
\DeclareMathOperator{\cl}{cl}
\DeclareMathOperator{\lin}{lin}

\def\C{\mathbb{C}}
\def\R{\mathbb{R}}

\def\N{\mathbb{N}}
\def\Z{\mathbb{Z}}

\def\l1{\ell^1}

\DeclareMathOperator{\Aw}{\mathcal{A}_\omega}

\newcommand{\vertiii}[1]{{\left\vert\kern-0.25ex\left\vert\kern-0.25ex\left\vert #1 
    \right\vert\kern-0.25ex\right\vert\kern-0.25ex\right\vert}} 

\title[Arens regularity of weighted convolution algebras that $ \cdots $ ]{Arens regularity of weighted convolution algebras that arise from totally ordered semilattices}

\bibliography{thesis_refs}

\begin{document}

\begin{abstract}
    We study Arens regularity of weighted semigroup convolution algebras for the specific case of totally ordered semilattices. This paper is a natural continuation of that of \cite{DalesStrauss}, where they studied the unweighted case. We also study some other properties of these algebras, such as the existence of Banach-algebra preduals.
\end{abstract}

\author{M. Eugenia Celorrio}

\renewcommand{\thefootnote}{\fnsymbol{footnote}} 
\renewcommand{\thefootnote}{\arabic{footnote}}

\keywords{Arens regularity, Arens products, weighted convolution algebra, totally ordered semilattice, dual Banach algebra, Banach algebra}

\maketitle

\section{Introduction}

\subsection{Setting the scene}
Arens regularity of Banach algebras continues to be a topic that draws the attention of Banach algebra specialists 70 years after Arens introduced the first and second Arens products in \cite{Arens-2}. In that same year, Arens proved in \cite{Arens-1} that the 
algebra $\ell^1(\N)$ with pointwise product is not Arens regular, which again seeded the start of a very prolific line of research on which we focus on this paper, the study of semigroup algebras.

Given a Banach algebra $A$, its product can be naturally extended to its bidual $A''$ in a way that $A$ is a closed subalgebra of $A''$. This extension can be done in two completely symmetrical ways. From these, the first ($\Box$) and second ($\Diamond$) Arens products arise. The biggest set on which these two products can agree is the whole bidual, while the smallest set on which they have to be equal is $A$. This sparked the definition of Arens regularity and later on strong Arens irregularity in \cite{DalesLau}, and the study of these two characteristics under different circumstances. In particular, Arens regularity of semigroup algebras has been broadly studied, for example see \cite{DalesLauStrauss}, \cite{DalesStrauss} where they studied semigroup algebras without a weight, while in \cite{Dedania}, \cite{filalitopcent} and \cite{Daws-Connes-Amenability} weighted convolution algebras are the main focus of the publications. In this paper we focus on the study of the weighted case. One of the earliest results about Arens regularity of semigroup algebras can be seen in \cite{CrawYoung}, where in 1974, Craw and Young studied Arens regularity of weighted semigroup algebras, focusing on cancellative semigroups.

\subsection{Overview of the paper}

We start the paper by giving, in Section 2, the relevant definitions of Arens products, Arens regularity, weighted semigroup algebras as well as any other background result that will be useful during the paper. We give an alternative proof of the necessary part of the main result \cite[Theorem 1]{CrawYoung} that reduces significantly the length of their proof. In this section we have also added some results that relate to general semigroups, and not the specific family that will be the focus of the paper.

We continue by introducing, in Section 3 the family of \textit{semilattices} that will be the main focus of the paper. The restriction to cancellative or weakly cancellative semigroups when studying Arens regularity of weighted semigroup algebras is very common. However, the family of semilattices introduced in this paper includes examples that are not necessarily weakly cancellative. So, our study includes a more general set-up in the study of Arens regularity. The inspiration to consider this specific family of semigroups came from \cite{DalesStrauss}, where they focus on studying Arens regularity of the unweighted case. Their main result is a characterization of the property of strong Arens irregularity in terms of some properties of $S$. However, we shall see in Example \ref{ex: Z wedge not SAI-AR} and Example \ref{ex:QnotSAI} that this characterization cannot be translated to weighted semigroup algebras. We provide a sufficient condition for a weighted semigroup algebra not to be strongly Arens irregular in Proposition \ref{prop: S wedge not SAI}. In contrast with strong Arens irregularity, in Theorem \ref{theo:Arens regularity when Limit of omega is infinity} we provide a characterization of Arens regularity of the weighted semigroup algebra, that depends solely on the properties of the weight. Banach algebras that have the extra property of being dual Banach algebras (see the definition in \S2) can be Arens regular under certain conditions. Hence it is natural to study this property (see Proposition \ref{theo:omegaboundeddual}, Proposition \ref{theo:dualSomegaunbounded} and \ref{prop:dualomegaboundedSwedge}). In the general study of dual Banach algebras, it is also interesting to see when a Banach-algebra predual is unique (in the sense that will be defined in Section 2). Thus, we have some results that talk about this, namely Proposition \ref{prop:S abelian idempotent. E_w=L(character space) implies unique predual} and Proposition \ref{prop: N unique predual}

In Section 4, we focus on a specific semilattice, $S=(\N, \wedge)$. This semilattice has attracted specialists for many years for the interesting properties that it has, for example in \cite{choi_2013} the focus is on the AMNM property. For us it is specially interesting because it is not weakly cancellative. In \cite{DalesLauStrauss} Arens regularity of the unweighted case of this semilattice was studied. In Proposition \ref{cor:DTCsetNminsemigroupalg} we obtain analogous results. Another reason why we focus on this family of weighted semigroup algebras is because, as it was seen in \cite{DalesLoy}, their Gel'fand transforms are weighted bounded variation algebras previously studied in \cite{Feinstein} and in \cite{White}.

\section{Preliminaries and general results}

\subsection{Banach Algebras preliminaries} 
In the following, we shall give a brief review of background, notation, and terminology that will be relevant to this paper. We shall follow mainly the conventions of \cite{DalesUlger}.

Given a Banach space $E$, the dual is denoted by $E'$. We start by introducing the two Arens products, $\Box$ and $\Diamond$ on the bidual of a Banach algebra $A$. We assume familiarity with the dual module, and we recall that the dual module operations are written as $a\cdot \lambda$ and $\lambda \cdot a$, for $a\in A$ and $\lambda \in A'$. For $\lambda \in A'$ and $M\in A''$, we can define $\lambda \cdot M$ and $M\cdot \lambda$ in $A'$ in the following way
    \begin{align*}
      \langle a,\lambda\cdot M \rangle=\langle M, a\cdot \lambda \rangle, \quad  \langle a,M\cdot \lambda \rangle=\langle M, \lambda\cdot a \rangle.  
    \end{align*}If $M\in A$ these new definitions agree with the original module operations. For $M,N \in A''$, we define $M\Box N$ and $M\Diamond N $ in $A''$ by
    \begin{align*}
        \langle M\Box N, \lambda \rangle=\langle M,N\cdot \lambda \rangle, \quad \langle M\Diamond N, \lambda\rangle=\langle N,\lambda \cdot M\rangle \quad (\lambda \in A'),
    \end{align*}
When both $M$ and $N$ are in $A$, then $M\Box N= M\Diamond N= MN$ (the usual product in $A$). It can be seen that this definition is the same as the following:

Let $M,N\in A''$, and take $(a_{\alpha}), (b_{\beta})$ nets in $A$ such that $\lim\limits_{\alpha} a_{\alpha}=M$ and $\lim\limits_{\beta} b_{\beta}=N$ in the weak-$*$ topology. Then the two Arens products are 
    \begin{align*} 
        M\Box N=\lim\limits_{\alpha} \lim\limits_{\beta} a_{\alpha}b_{\beta},\quad M\Diamond N=\lim\limits_{\beta} \lim\limits_{\alpha} a_{\alpha}b_{\beta},
    \end{align*}
where the limits are again in the weak-$*$ topology on $A''$.

We assume familiarity with the bidual module, and we recall that, for $a\in A$ and $M\in A''$, we write by $a\cdot M$ and $M\cdot a$ the module operations. The two Arens products defined as above are one sided $\sigma(A'', A')-$continuous, each on a different side. In particular, if $(s_\alpha)$ is \ a net in $A$ such that $(s_\alpha)\rightarrow M$ in the weak-$*$ topology and $N\in A''$, then $s_\alpha \cdot N\rightarrow M\Box N$ and $N\cdot s_\alpha\rightarrow N\Diamond M$, both limits in the weak-$*$ topology again. Notice that we are identifying $A$ with the canonical image in its bidual.
Throughout the rest of this paper, unless specified otherwise, whenever we talk about the bidual of a Banach algebra we are implicitly talking about $(A'',\Box)$.

\begin{deff}
    Let $A$ be a Banach algebra. The \textit{left topological centre} of $A''$ is
        \begin{align*}
            \mathfrak{Z}^{(\ell)}(A'')=\{ M\in A'' : M\Box N= M\Diamond N \quad(N\in A'')\}.
        \end{align*}
    Similarly, we define the \textit{right topological centre} of $A''$ as
        \begin{align*}
            \mathfrak{Z}^{(r)}(A'')=\{ M\in A'' : N\Box M= N\Diamond M \quad (N\in A'')\}.
        \end{align*}
    These two topological centres might be different. However, in the case where $A$ is commutative, the right and the left topological centre are the same and we can speak about the \textit{topological centre} of $A''$, which is
        \begin{align*}
            \mathfrak{Z}(A'')=\{ M\in A'' : M\Box N= M\Diamond N \quad (N\in A'')\}.
        \end{align*}
\end{deff}
Since, for $a\in A$ and $N\in A''$, we have $N\Box a= N\Diamond a= N\cdot a$ and $a\Box N= a\Diamond N= a \cdot N$, we always have $A\subset \mathfrak{Z}^{(\ell)}(A'') \subset A''$ (also $A\subset \mathfrak{Z}^{(r)}(A'') \subset A''$). It might be that the left or right topological centres are neither $A$ nor $A''$, which leads to the following two definitions:
\begin{deff}
    Let $A$ be a Banach algebra. We say $A$ is  \textit{Arens regular} when     
        $$
            \mathfrak{Z}^{(\ell)}(A'')=\mathfrak{Z}^{(r)}(A'')=A''.
        $$
    A Banach algebra $A$ is said to be \textit{strongly Arens irregular} if
        $$
            \mathfrak{Z}^{(\ell)}(A'')=\mathfrak{Z}^{(r)}(A'')=A.
        $$
    
    In the special case when $A$ is a commutative Banach algebra, $A$ is \textit{Arens regular} if and only if $\mathfrak{Z}(A'')=A''$ and \textit{strongly Arens irregular} if and only if $\mathfrak{Z}(A'')=A$.
\end{deff}
Thus, a commutative Banach algebra is Arens regular if and only if $(A'', \Box)$ is commutative.

Let $A$ be a Banach algebra. Given $\lambda \in A'$ we say that $\lambda $ is \textit{weakly almost periodic} if 
    \begin{align*}
        R_\lambda: A&\rightarrow A'\\
        a&\mapsto a\cdot \lambda
    \end{align*}
is weakly compact. Following \cite{DalesUlger} we write $\WAP(A)$ for the space of weakly almost periodic functionals on $A$.

It is standard that $A$ is Arens regular if and only if $\WAP(A)=A'$. See \cite{Lust-piquard}.   

We shall introduce now the definition of Banach-algebra predual, following \cite{DalesUlger}. We shall see below how the notion of dual Banach algebra is closely tied to that of Arens regularity.

\begin{deff}\label{def:Banach preduals}
    Let $E$ be a Banach space. Then a closed subspace $F$ of $E'$ is a \textit{concrete predual of $E$} if the map $T_F: E \rightarrow F'$ defined by 
        $$
         (T_Fx)(\lambda)=\langle x, \lambda \rangle _{E,E'} \quad (x\in E, \lambda \in F)
        $$
    is a linear homeomorphism.

    Let $A$ be a Banach algebra. A \textit{Banach-algebra predual} for $A$ is a closed linear subspace $F$ of $A'$ that is a concrete predual of $A$ and an $A$-bimodule. We say that $A$ is a \textit{dual Banach algebra} if it has a Banach-algebra predual. A Banach-algebra predual is \textit{unique} if it is the only closed submodule of $A'$ with respect to which $A$ is a dual Banach algebra.
\end{deff}

It is standard that for a Banach space $E$ and a concrete predual $F\subset E'$, we can write $E''$ as the following direct sum: $E''=E\oplus F^\perp$, where we are identifying $E$ with its canonical image on $E''$ and $F^\perp\subset E''$ is the usual orthogonal complement.

The following proposition can be found in \cite[Proposition 1.3.25]{DalesUlger}.

\begin{prop}\label{prop: concrete preduals subset}
    Let $E$ be a Banach space, and suppose that $F$ and $G$ are concrete preduals of $E$ such that $F\subset G$. Then $F=G$. \hfill $\Box$
\end{prop}

\subsection{Weighted semigroup algebras}

We shall introduce here some terminology, definitions and background results. We shall also introduce some generic results that will be useful in later sections, such as Proposition \ref{prop:character space of weighted semigroup algebras (idempotents)}.

\begin{deff}
    Let $S$ be a semigroup.  We say that $S$ is \textit{right cancellative} (respectively, \textit{left cancellative}) if, for all $a,s,t \in S$, $sa=ta$ (respectively, $as=at$) implies that $s=t$. When $S$ is both right and left cancellative we call it \textit{cancellative}. 
    
    We say $S$ is \textit{weakly right cancellative} (respectively, \textit{weakly left cancellative}) if, for all $s,t \in S$, the set $\{ u\in S: us=t \}$ (respectively, $\{ u\in S: su=t \}$)  is finite. If $S$ is both right and  left weakly cancellative, we say that $S$ is\textit{ weakly cancellative}.
    
    An element $p \in S$ is an \textit{idempotent} if $p^2=p$. We say that $S$ is an \textit{idempotent semigroup} if every element of $S$ is idempotent. We say that a semigroup $S$ is \textit{separating} if $s=t$ whenever $s,t \in S$ are such that $st=s^2=t^2$.
\end{deff}  
Whenever $S$ is cancellative or idempotent, then $S$ is separating.

We recall that a semilattice is a partially ordered set such that every nonempty finite subset has a greatest lower bound. In our context we have the following:
\begin{deff}\label{def: semilattices}
    Let $S$ be an abelian idempotent semigroup. Then we can define a partial order $\leq$ in $S$ by setting
        \begin{align*}
            s\leq t \iff st=s \quad (s,t \in S).
        \end{align*}
    For every $s,t \in S$, it can be seen that $st$ is a greatest lower bound for $\{s,t\}$. Hence, $(S, \leq)$ is a semilattice. Symmetrically, when we have a semilattice $(S,\leq)$, we can define a semigroup operation by setting $st$ as the greatest lower bound of $\{s,t\}$ $(s,t \in S)$. Hence, from now on, we shall say that $S$ is a \textit{semilattice} when $S$ is an abelian idempotent semigroup.
\end{deff}

In order to talk about weighted semigroup algebras, we first need to introduce the notion of a weight on a semigroup.

 \begin{deff}
     Let $S$ be a semigroup. A \textit{weight on }$S$ is a function $\omega:S \longrightarrow (0,\infty)$ such that it is submultiplicative, in the sense that 
        $$
            \omega(st)\leq \omega(s)\omega(t) \quad (s,t \in S).
        $$
 \end{deff}
When $S$ is an idempotent semigroup, $\omega:S \longrightarrow (0,\infty)$ is a weight if and only if $\omega(s)\geq 1$, $(s\in S)$.

Let $S$ be a semigroup. We shall denote by $\delta_s$ the characteristic function of an element $s\in S$. Given a weight $\omega$ on $S$, we define $\tilde{\delta}_s\coloneqq \delta_s/\omega(s).$

\begin{deff}
    Let $S$ be a semigroup and let $\omega$ be a weight on $S$. Then we define the \textit{weighted semigroup algebra} of $S$ as the Banach space
        \begin{align*}
            \mathcal{A}_\omega:=\ell^1(S,\omega)= \left\lbrace \alpha=\sum\limits_{s\in S} \alpha(s) \delta_s: \|\alpha\|_\omega \coloneqq \sum\limits_{s\in S} |\alpha(s)|\omega(s) <\infty \right\rbrace,
        \end{align*}
    where $\alpha(s)\in \mathbb{C}$ ($s\in S$), together with the convolution multiplication that can be determined by the partial definition given below, via continuous bilinear extension:
        $$
            \delta_s\star \delta_t = \delta_{st} \quad (s,t \in S).
        $$
\end{deff}

For $\omega \equiv 1$, this is the usual convolution algebra. We shall refer to this specific situation as the \textit{unweighted case}.

\begin{rem}\label{rem: linf dual Ew predual}
      The dual of $\mathcal{A}_\omega$ is 
        \begin{align*}
            \mathcal{A}_\omega':=\ell^\infty(S,1/\omega)= \left\lbrace \lambda \in \mathbb{C}^S: \sup\{|\lambda(s)|/\omega(s) : s\in S\}<\infty \right\rbrace,
        \end{align*}
    with the norm denoted by $\|\cdot\|'_\omega$ so that
        \begin{align*}
            \|\lambda\|'_\omega : = \sup\{|\lambda(s)|/\omega(s) : s\in S\} \quad (\lambda \in \ell^\infty(S,1/\omega)).
        \end{align*}
    The space 
        $$
            E_\omega:=c_0(S,1/\omega)=\left\lbrace \lambda \in \mathbb{C}^S: \lim\limits_{s\rightarrow \infty} |\lambda(s)|/\omega(s)=0\right\rbrace
        $$
    is a concrete predual of $\mathcal{A}_\omega$, as in Definition \ref{def:Banach preduals}.
\end{rem}
In the following section we shall see that $E_\omega$ is not necessarily a Banach-algebra predual for $\Aw$, and we shall study in which cases it is.

The following lemma is straightforward, and we omit the proof. We add it here to facilitate the reading of the paper: 

\begin{lem}\label{lem: theta_omega} 
    Let $S$ be a semigroup and let $\omega$ be a weight on $S$. Let 
    \begin{align*}
        \theta_\omega: \alpha\mapsto \alpha/\omega, \text{ } \ell^1(S)\longrightarrow  \Aw,
    \end{align*}
    where $$
    \alpha/\omega\coloneqq \sum\limits_{s\in S} \alpha(s) \delta_s/\omega(s), \quad (\alpha=\sum\limits_{s\in S} \alpha(s) \delta_s).
    $$
    Then $\theta_\omega$ is an isometric isomorphism of Banach spaces. Note that, for every $ s\in S$, $\theta_\omega(\delta_s)=\tilde{\delta}_s$ . \hfill $\Box$
\end{lem}

Given a semigroup $S$, we denote by $\varPhi_S$ the space of semi-characters on $S$. Recall that a semicharacter on $S$ is a map \mbox{$\theta: S \longrightarrow \bar{\mathbb{D}}$} such that $\theta\neq 0$ and
    $$
        \theta (st)=\theta(s)\theta(t) \quad (s,t \in S).
    $$ 
    For more information on semicharacter of semigroups in this setting we recommend \cite{DalesLauStrauss}. Given a Banach algebra $A$, we denote the character space of $A$ by $\varPhi_A$. Let $S$ be a semigroup and $\omega$ a weight on $S$. The character space of the weighted semigroup algebra $\Aw$ is denoted by $\varPhi_\omega$. 

 For a weight bounded below, we have that $\varPhi_S\subset \varPhi_\omega$ as a subset. It is well known, see e.g. \cite[\S4]{DalesLauStrauss}, that for $\omega \equiv 1$ this inclusion is an equality. For general weights and semigroups the inclusion could be proper. For example, take $S=\mathbb{Z}$ with addition and $\omega(n)=\mathrm{e}^{|n|}$, in which case the semicharacter space, $\varPhi_S$, is the unit circle while the character space, $\varPhi_\omega$, is an annulus with inner radious $r=e^{-1}$ and outer radious $R=e$. The details are left to the reader. However, as we shall see in the result below, for the semigroups in this paper this problem does not arise.

\begin{prop}\label{prop:character space of weighted semigroup algebras (idempotents)} 
    Let $S$ be a semilattice, and let $\varPhi_S$ be the semicharacter space of $S$. Let $\omega: S \rightarrow [1,\infty)$ be a weight on $S$. Then $\varPhi_\omega=\varPhi_S$.
\end{prop}
\begin{proof}
    We know that $\varPhi_S \subset \varPhi_\omega$. Now let $\varphi$ be a character on $\Aw$ and define $\theta_\varphi(s)=\varphi(\delta_s)$. Since $\varphi$ is a character, 
        $$
            \theta_\varphi(s)=\varphi(\delta_s)=\varphi(\delta_s\star \delta_s) = \varphi(\delta_s)\varphi(\delta_s) \quad(s\in S),
        $$
    and so $\theta_\varphi(s)\in \{0,1\}$. Also
        $$
            \theta_\varphi(st)=\varphi(\delta_{st}) = \varphi(\delta_s \star \delta_t)=\varphi(\delta_s)\varphi(\delta_t)=\theta_\varphi(s)\theta_\varphi(t)\quad(s,t \in S).
        $$ 
    Hence $\theta_\varphi$ is a semi-character on $S$ and so we can identify $\varPhi_\omega$ with $\varPhi_S$.
\end{proof}

Following \cite{DalesLau} we introduce the following definition.

\begin{deff}\label{def:limites}
    Let $S$ be an infinite set. Let $f:S\longrightarrow \mathbb{R}$. Given $C \in \mathbb{R}$, we write
        $$
          \Lim\limits_{s\rightarrow \infty} f(s) =C
        $$
    if, for each $\varepsilon>0$, there is a finite set $F$ of $S$ such that
        $$
            |f(s)-C|<\varepsilon \quad (s\in S\setminus F).
        $$
    We write
        $$
            \Lim\limits_{s\rightarrow \infty} f(s)=\infty
        $$
    if, for each $M>0$, there is a finite set $F$ of $S$ such that 
        $$
            f(s)>M \quad (s \in S \setminus F).
        $$
    We write
        $$
            \Liminf\limits_{s\rightarrow \infty} f(s) <\infty 
        $$
    if and only if it is not true that $\Lim\limits_{s\rightarrow \infty} f(s) = \infty$, i.e., there exists $M>0$ such that the set $\{ s\in S: f(s)<M\}$ is infinite. We write $\Liminf f$ for the infimum of these constants $M$.
\end{deff}

We denote by $\beta S$ the Stone-\v{C}ech compactification of $S$, where we have given $S$ the discrete topology. We denote by $S^*$ the growth of $S$, which is defined to be $\beta S \setminus S.$

We proceed to introduce two definitions that are key in the study of Arens regularity of Banach algebras.

\begin{deff}\label{def:Omega} 
    Given a weight $\omega$ on $S$,we define $\Omega$ on $S\times S$ in the following way:
        $$
            \Omega(s,t)=\frac{\omega(st)}{\omega(s)\omega(t)} \quad (s,t \in S).
        $$
    Given a function $f:S\times S\rightarrow \R$, we say that $f$ \textit{$0-$clusters on $S \times S$} if, for any $(x_n), (y_m)$ sequences, each consisting of distinct elements of $S$, then
        \begin{align*}
            \lim_{n\rightarrow \infty} \lim_{m\rightarrow \infty} f(x_n, y_m)=\lim_{m\rightarrow \infty} \lim_{n\rightarrow \infty} f(x_n, y_m)=0
        \end{align*}
    whenever both repeated limits exist.
\end{deff}
For more on repeated limit conditions, we recommend \cite[\S 3]{DalesLau}.

In 1974 Craw and Young already studied Arens regularity of weighted semigroup algebras in \cite{CrawYoung}. We provide a new proof of their main theorem (\cite[Theorem 1]{CrawYoung}) below. The proof has been modified to match the terminology used here, which resulted in a simpler version of the necessary part. A similar observation was made in \cite{Daws-Connes-Amenability}, however, we provide full details here.

\begin{theo}[Craw-Young,1974]\label{theo:crawyoung}
    Let $S$ be an infinite semigroup and $\omega$ a weight on $S$. Then:
    \begin{enumerate}
        \item[(a)] If $\Omega$ $0-$clusters on $S\times S$, then $\Aw$ is Arens regular. 
        \item[(b)] If $S$ is a weakly cancellative semigroup, then the Arens regularity of $\Aw$ implies that $\Omega$ $0-$clusters.
    \end{enumerate}
\end{theo}
\begin{proof}
     $(a)$ Let $\lambda \in \Aw'$. Then $\lambda\in \WAP(\Aw)$ if and only if
        \begin{align*}
            (i,j)&\mapsto \langle R_\lambda (\tilde{\delta_i}), \tilde{\delta}_j \rangle , \quad t_\omega: S\times S \rightarrow \mathbb{C}
        \end{align*}
    clusters. This follows from Grothendieck's double limit criterion and can be seen in \cite{Daws-Connes-Amenability}. This translates to the fact that
        $$
            t_\omega(i,j)=\frac{\langle \delta_i \cdot \lambda ,\delta_j\rangle}{\omega(i)\omega(j)}=\Omega(i,j)\langle \lambda, \delta_i\star \delta_j\rangle
        $$
    clusters.

    \noindent In particular, when $\Omega$ $0-$clusters, $t_\omega$ also $0$-clusters. Thus, for $\lambda \in \Aw'$, \mbox{$\lambda \in \WAP(\Aw)$}. Hence, $\Aw$ is Arens regular.

    $(b)$ Suppose that $S$ is weakly cancellative. The following argument is very similar to the one followed in \cite[Theorem 1]{CrawYoung} but we provide full details for completeness. Suppose that $\Omega$ does not $0$-clusters, i.e. there exist $(s_n)$ and $(t_m)$ sequences, each consisting of distinct elements of $S$, such that the repeated limits exists and at least one of them is not equal to $0$. We may suppose $\lim\limits_{n\rightarrow \infty} \lim\limits_{m\rightarrow \infty} \Omega(s_n, t_m)\neq 0$. Take $\varepsilon>0$ such that
        $$
	        \lim\limits_{n\rightarrow \infty} \lim\limits_{m\rightarrow \infty} \Omega(s_n, t_m) > \varepsilon.
        $$
    We shall see that there are two elements of the bidual such that the two Arens products are different.
    We may suppose that 
    $$
        \lim\limits_{m\rightarrow \infty} \Omega(s_n, t_m)>\varepsilon>0 \quad (n\in \N).
    $$
    Let us choose two subsequences $(s'_n)$ and $(t'_m)$ of $(s_n)$ and $(t_m)$, respectively, such that $\Omega(s'_n,t'_m)>\varepsilon$ for $n\leq m$. Indeed, take $s'_1=s_1$ and take $t'_1$ to be the first element $t_m$ such that $\Omega(s_1,t_m)>\varepsilon.$ Let us suppose that we already  and consider the set $F$ of elements $u\in S$ such that $ut'_l=s'_it'_j$ for at least one triple $l,i,j$ with $1\leq j<k, 1\leq i,l\leq k$. Since $S$ is weakly cancellative then $F$ is finite. Hence we can choose as $s'_k$ the first element $s_n$ such that $s_n\notin F$ with $n> k$ ($u\in F$). Following the same line of reasoning, we consider the set $E$ of elements $u\in S$ such that $s'_l u=s'_i t'_j$ for at least one triple $l,i,j$ with $1\leq j<k, 1\leq i,l\leq k$, which is also finite. So we can choose $t'_k$ as the first element $t_m$ with $m>k$ such that 
        \begin{align*}
            \Omega(s'_i,t_m)>\varepsilon \quad (1\leq i\leq k),
           \quad t_m\notin E.
        \end{align*}
    These subsequences are such that $\Omega(s'_n, t'_m) >\varepsilon$ for $n\leq m$ and such that the elements $s'_n t'_m$ are distinct for $m,n \in \N$. Let $\alpha_n$ and $\beta_m$ be the normalized point masses at $s'_n$ and $t'_m$ respectively and $\chi\in \ell^\infty (S)$ the characteristic function of the set $\{s'_n t'_m: n\leq m\}$. Then 
        \begin{align*}
            \langle \alpha_n \star \beta_m, \omega \chi\rangle  =\Omega(s'_n,t'_m)>\varepsilon &\quad (n\leq m),\\
            \langle \alpha_n \star \beta_m, \omega \chi\rangle =0 &\quad (n>m).
        \end{align*}
    Let $M,N \in \Aw ''$ be $\sigma(\Aw'',\Aw')$-accumulation points of $(\alpha_n)$ and $(\beta_n)$ respectively. By construction, $\langle M\Box N, \omega \chi\rangle\geq \varepsilon$ and  $\langle M\Diamond N, \omega \chi\rangle=0$. Thus, $\Aw$ is not Arens regular, as desired.

\end{proof}

The following result and full details of the proof can be found in \cite[Proposition 3.1]{DalesLau}. We write it here for completeness.

\begin{prop}\label{prop:Omegaboxanddiamond}
    Let $S$, $T$ be non-empty sets. Take \mbox{$f: S\times T\rightarrow \C$}. Suppose that $(s_\alpha) $ and $(t_\beta)$ are nets in $S$ and $T$, respectively, such that $a=\lim_\alpha \lim_\beta f(s_\alpha,t_\beta)$ and $b=\lim_\beta \lim_\alpha f(s_\alpha,t_\beta)$ both exist. Then there are subsequences $(s_{\alpha_m})$ and $(t_{\beta_m})$ of the nets $(s_\alpha) $ and $(t_\beta)$, respectively, such that $a=\lim_m \lim_n f(s_{\alpha_m},t_{\beta_n})$ and $b=\lim_\beta \lim_\alpha f(s_{\alpha_m},t_{\beta_n})$.
\hfill $\Box$
\end{prop} 
By subsequence we mean that $\alpha_1 \preceq \alpha_2 \preceq \alpha_3 \cdots$.

Consider $\Omega: S\times S\longrightarrow \mathbb{R}$ defined as in Definition \ref{def:Omega}. Recall that, for $s,t\in S$, we have $0\leq \Omega(s,t) \leq 1$.  Let $u,v\in \beta S$. Then there are nets $(s_\alpha), (t_\beta)$ in $S$ such that $u=\lim_\alpha s_\alpha$, $v=\lim_\beta t_\beta$. Recall from \cite[\S 4]{Dedania} that we can define
    \begin{align*}
        \Omega_\Box (u,v)=\lim_\alpha\lim_\beta \Omega(s_\alpha,t_\beta),\quad
        \Omega_\Diamond (u,v)=\lim_\beta\lim_\alpha \Omega(s_\alpha,t_\beta).
    \end{align*}
We can see that if we apply Proposition \ref{prop:Omegaboxanddiamond} to $\Omega_\Box$ and $\Omega_\Diamond$ we can define them in terms of sequences instead of nets.

\begin{rem}\label{rem:Box Diamond product general weighted algebra}
    Let $S$ be a semigroup and let $\omega$ be a weight on $S$.Let $u,v \in S^*$ and consider $\theta_\omega$ as in Lemma \ref{lem: theta_omega}. Then
    \begin{align*}
        \langle\theta_\omega''(\delta_u)\Box_\omega \theta_\omega''(\delta_v),\omega\lambda\rangle=\Omega_\Box(u,v)\langle\delta_u\Box \delta_v,\lambda\rangle,\\ \langle\theta_\omega''(\delta_u)\Diamond_\omega \theta_\omega''(\delta_v),\omega\lambda\rangle=\Omega_\Diamond(u,v)\langle\delta_u \Diamond \delta_v,\lambda\rangle.
    \end{align*}
    This follows from a simple calculation, using nets $(s_\alpha)$ and $(t_\beta)$ in $S$ such that $u=\lim_\alpha s_\alpha$ and $v=\lim_\beta t_\beta$.
\end{rem}

\section{Results for totally ordered semilattices}

Let $(S,\leq)$ be a semilattice, as described in Section 2. Suppose also, that the order $\leq$ is a total order, meaning that for any two elements $s,t\in S$, it is always true that either $s\leq t$ or $t\leq s$. We shall refer to $(S,\leq)$ with these characteristics as a \textit{totally ordered semilattice}. Totally ordered semilattices are the object of study of this section. Conversely, any totally ordered set $S$ becomes a semigroup if we take $st=\min\{s,t\}$; this semigroup is a semilattice and the partial order defined in Definition \ref{def: semilattices} coincides with the original order on $S$.

Some preliminaries about semigroups of the form of $S$ are in \cite{Ross}.

\begin{rem}
Notice that the natural numbers with the usual order $(\N, \leq)$ belong to this family. This semilattice arises from the semigroup $\N$ with the minimum operation, which is not even weakly cancellative. Thus, the family studied comprises a wide variety of examples with characteristics that differ from those usually considered in the literature focused on weighted semigroup algebras. We shall focus on this example in the following section.
\end{rem}

\subsection{Arens regularity}

We shall assume from now on that there exists an embedding from $S$ into some infinite, semigroup $T$ that contains $S$ as a subsemigroup, with the following aditional characteristics:
    \begin{itemize}
        \item $T$ must be a totally ordered set that preserves the order in $S$;
        \item $T$ has a minimum and a maximum, which we shall call $0$ and $\infty$, respectively;
        \item $T$ is complete in the sense that every non-empty subset of $T$ has a supremum and an infimum;
        \item We consider the interval topology on $T$, in which case $T$ is a compact topological semigroup.
    \end{itemize}

Note that in this case every strictly increasing, respectively strictly decreasing, net in $S$ converges to its supremum, respectively infimum. Note that the notion of limit here used is different from the previously introduced $\Lim$.

In the following remark we shall see that given an infinite totally ordered semilattice $S$, we can always find such a $T$. This is well-known to specialists, but we add it to make it more accessible for the reader.

\begin{rem}\label{rem: S can be embedded in T}
 Since $(S,\wedge)$ is a semilattice, we know that $\varPhi_S$ separates the points of $S$. Let $\varSigma$ be a subset of $\varPhi_S$ that separates the points of $S$, and let $\kappa= |\varSigma|$. Since every character of $S$ maps into $\{0,1\}$, $(S, \wedge)$ can be embedded as a semigroup in  $C:=(\{0,1\}^\kappa,\wedge)$. Let $T$ be the closure of $S$ in $C$. In this case, $T$ is a complete totally ordered lattice which is compact in its interval topology, as desired. For details, included the definition of the order of $T$ (which is induced by the order on $S$), see \cite[\S2]{DalesStrauss}.
\end{rem}

Let $S$ be a totally ordered semilattice and $T$ as above. Let $U$ be a subset of $S$, we write as $\cl_T U$ and $\cl_{\beta S} U$ the closures of $U$ in $T$ and in $\beta S$, respectively. The continuous extension of the inclusion map of $S$ into $T$ is denoted by
    $$
        \pi:\beta S \rightarrow T.
    $$
Thus $\pi(\beta S)=\cl_T S$. For $t\in \cl_T S$, we shall write $F_t$ for the fibre $\{x \in \beta S: \pi(x)=t\}$ and $F^*_t=F_t \cap S^*$ (recall that $S^*$ is the growth of $S$).
Thus we have that 
    $$
        F^*_t= F_t \quad (t\in T\setminus S) \quad\text{and} \quad F^*_t=F_t\setminus \{t\} \quad (t \in S).
    $$

We shall denote by $E$ the set of accumulation points of $S$ in $T$. We have then that $E\neq \emptyset$. For $t\in T$, $F_t^*$ is a closed,
compact subset of $\beta S$ and $F_t^* \neq \emptyset$ if and only if $t\in E$.

Consider $\Omega$ as in Definition \ref{def:Omega}. When $S$ is a totally ordered semilattice, $\Omega(s,t)=1/\omega(t)$ $(s\leq t)$. What is more, we can see that $\Omega$ $0-$clusters if and only if $\Lim\limits_{s\rightarrow \infty} \omega(s)=\infty.$ Indeed $\Lim\limits_{s\rightarrow \infty} \omega(s)=\infty$ implies that $\Omega$ $0-$clusters, since, for every $\varepsilon>0$, the set of elements $s\in S$ such that $\frac{1}{\omega(s)}>\varepsilon$ is finite. 

\noindent Now suppose that $\Liminf\limits_{s\rightarrow \infty} \omega(s) <\infty $. Then there exists $M>0$ such that the set $\{ s\in S: \omega(s)<M\}$ is infinite. Hence, we can then take two sequences $(s_m), (t_n)$ of distinct elements belonging to that set and so $\Omega(s_m,t_n)\geq 1/M^2 >0$. Since $\Omega$ is always bounded above by $1$, we can choose subsequences $(s_{m_k})$ and $(t_{n_j})$ such that the repeated limit $\lim_m \lim_n \Omega(s_{m_k},t_{n_j})$ exists. We can repeat this to find subsequences of $(s_{m_k})$ and $(t_{n_j})$ such that the repeated limit with reverse indexes exist. Since both these limits are greater than $1/M^2 >0$, we have that $\Omega$ does not $0-$cluster. 

\begin{theo}\label{theo:Arens regularity when Limit of omega is infinity}
    Let $(S, \wedge)$ be an infinite, totally ordered semilattice. Let $\omega$ a weight on $S$. Then the following conditions are equivalent:
        \begin{enumerate}
            \item[(a)] the algebra $\Aw$ is Arens regular;
            \item[(b)] $\Lim\limits_{s\rightarrow \infty} \omega(s)=\infty$;
            \item[(c)] $M\Box N=M\Diamond N=0$ $(M,N \in E_\omega^\perp).$
        \end{enumerate}
\end{theo}
\begin{proof}
    $(a) \Rightarrow (b)$ Suppose that $\Liminf_{s\rightarrow \infty} \omega(s) <\infty$, and let $M> \Liminf\omega$ . Let \mbox{$U\coloneqq \{ s\in S: \omega(s)<M\}$}. Take $t\in E\cap \cl_{\beta S} U$. By \cite[Lemma 2.4]{DalesStrauss}, $|F_t^{*}|=2^{\mathfrak{c}}$. Take $p\in F^*_t$. If $p\in F_t^* \cap \cl_{\beta S} (U \cap [0,t))$, there exists $q$ in $F_t^* \cap \cl_{\beta S} (U \cap [0,t))$ with $\delta_p \notin \lin \{ \delta_q \}$. Then 
        $$
            \delta_p\Box \delta_q = \delta_p \quad \text{and} \quad \delta_p\Diamond \delta_q = \delta_q.
        $$
    Note that this multiplications are in the bidual of the unweighted algebra. We shall look at what happens with the weighted case now.

    Let us consider the isometric isomorphism $\theta_\omega$ as in Lemma \ref{lem: theta_omega}. For $\lambda \in C (\beta S)$, we have
        \begin{align}\label{eq:Boxprod}
            \langle\theta_\omega(\delta_p)\Box_\omega \theta_\omega(\delta_q),\omega\lambda\rangle=\Omega_\Box(p,q)\langle \delta_p,\lambda\rangle,
        \end{align}
    And that
            \begin{equation}\label{eq:Diamondprod}
                \langle\theta_\omega(\delta_p)\Diamond_\omega \theta_\omega(\delta_q),\omega\lambda\rangle=\Omega_\Diamond(p,q)\langle \delta_q,\lambda\rangle.
            \end{equation}

    Take $(s_\alpha)$, $(t_\beta)$ nets in $S$ converging to $p$ and $q$, respectively. Observe that, since $1\leq \omega(s_\alpha)\leq M$ and $1\leq \omega(t_\beta)\leq M$ for all $\alpha, \beta$, we have that 
        $$
        	0<1/M^2\leq \Omega_{\Box} (p,q)\leq M \quad \text{and} \quad 0<1/M^2 \leq \Omega_\Diamond (p, q)\leq M.
        $$
    Hence the equations (\ref{eq:Boxprod}) and (\ref{eq:Diamondprod}) from above are equal for every $\lambda \in C(\beta S)$ if and only if $p= \frac{\Omega_\Diamond (p,q)}{\Omega_\Box (p,q)}q$. But that is not possible as $\delta_p\notin \lin\{\delta_q\}$. Thus, $\Aw$ is not Arens regular.

        $(b)\Rightarrow (c)$ follows from \cite[Theorem 6.3.23]{DalesUlger}.   

        $(c)\Rightarrow (a)$ Since $E_\omega$ is a concrete predual or $\Aw$, we know that $\Aw''=\Aw\oplus E_\omega^\perp$. Since $M\Box N=M\Diamond N=0$ for every $M,N \in E_\omega^\perp$, the result follows.
\end{proof}

The next question would be which conditions on $\omega$ ensure that $\Aw$ is strongly Arens irregular. In \cite{DalesStrauss}, they proved that the semigroup algebra $(\ell^1(S),\star)$ is strongly Arens irregular if and only if $\cl_T S$ is scattered. We shall see below that when we add a weight the situation is more complex. We start by considering the simplest case. When $S$ is a semigroup and $\omega$ is a bounded weight on $S$, then the inclusion map $\Aw \hookrightarrow \ell^1(S)$ is a Banach algebra isomorphism, and so $\Aw$ is strongly Arens irregular if and only if $\cl_T S$ is scattered.

However, we shall see below that if $\omega$ is not bounded we can have several different options.The next provides a sufficient condition for $\Aw$ to not be strongly Arens irregular.

\begin{prop}\label{prop: S wedge not SAI}
    Let $S$ as above and let $\omega$ be a weight on $S$. Suppose that for every $p\in F_\infty^*$ and every net $(s_\alpha)$ such that $s_\alpha \rightarrow p$, the set $\{\omega(s_\alpha)\}$ is unbounded. Then $\Aw$ is not strongly Arens irregular. 
\end{prop}
\begin{proof}
    Let $p\in F_\infty^*$. We shall see that $\delta_p\in \mathfrak{Z}(\Aw'')$. Take $q\in S^*$. If $q\in F_\infty^*$ then $\Omega_\Box(p,q)=\Omega_\Diamond(p,q)=0$ and so 
        $$
            \theta_\omega(\delta_p)\Box_\omega \theta_\omega(\delta_q)=0=\theta_\omega(\delta_p)\Diamond_\omega \theta_\omega(\delta_q).
        $$
    Take now $q \in S^*$ such that $q\notin F_\infty^*$ and let $(s_\alpha)$ a net in $S$ converging to $p$ and $(t_\beta)$ a net in $S$ converging to $q$. Since $\pi(q)<\pi(p)$, this implies that, deleting a finite number of elements if needed, 
    we can suppose that $t_\beta<s_\alpha$ for every $\alpha$, $\beta$. Hence we have that
        $$ 
            \Omega(t_\beta,s_\alpha)=\frac{\omega( t_\beta \wedge s_\alpha )}{\omega(t_\beta)\omega(s_\alpha)}=\frac{1}{\omega(s_\alpha)}
        $$
    and so $\Omega_\Box(p,q)=\Omega_\Diamond(p,q)=0$ which gives us again that
        $$
            \theta_\omega(\delta_p)\Box_\omega \theta_\omega(\delta_q)=0=\theta_\omega(\delta_p)\Diamond_\omega \theta_\omega(\delta_q).
        $$ 
   We conclude then that $\theta_\omega(\delta_p)\in \mathfrak{Z}(\Aw'')$ but $\theta_\omega(\delta_p)\notin \Aw$. Thus $\Aw$ is not strongly Arens irregular.
\end{proof}

This previous result, together with Theorem \ref{theo:Arens regularity when Limit of omega is infinity}, allow us to obtain plenty of weighted semigroup algebras that are neither Arens regular nor strongly Arens irregular. The following two examples portray two different semilattices, one of them is such that $\cl_T S$ is scattered and the other one is such that $\cl_T S$ is not scattered. In contrast with the unweighted case, we shall see that both of them are neither Arens regular nor strongly Arens irregular.

\begin{ex} \label{ex: Z wedge not SAI-AR}
    Let $S=\Z$, and $T=\{-\infty\}\cup \R \cup \{\infty\}$. Then, $\cl_T S= \{-\infty\}\cup \Z \cup \{\infty\}$, which is scattered.

    Consider $\omega$ a weight on $S$ such that $ \lim_{n\rightarrow \infty} \omega(n)=\infty $ and such that $\omega|(\Z \setminus \N)$ is bounded. Then $\Aw$ is neither Arens regular nor strongly Arens irregular. As an specific example, take $\omega(n)=n$, for $n\geq 1$, and $\omega(n)=1$, for $n<1.$
\end{ex}

\begin{ex}\label{ex:QnotSAI}
    Let $S= \mathbb{Q}^{+ \bullet}\coloneqq \{p\in \mathbb{Q}: p>0\}$. Consider a weight $\omega: \mathbb{Q} \rightarrow [1,\infty)$ such that $\omega(p)=1$ $(p \in [0,1] \cap S)$ and such that $\lim\limits_{p\rightarrow\infty} \omega(p)=\infty$. Then $\Aw$ is not Arens regular neither strongly Arens irregular.
    
    The fact that $\Aw$ is not strongly Arens irregular follows from Proposition \ref{prop: S wedge not SAI}. However, in this case we can find a concrete element $M\in  \mathfrak{Z}(\Aw'')$, but such that $M\notin \Aw$. Since $\lim\limits_{p\rightarrow\infty} \omega(p)=\infty$, there exists $(p_n)$ a strictly increasing sequence such that $\lim\limits_{n\rightarrow\infty}\omega(p_n)=\infty$. For clarity we will call that sequence $P$. Consider $u \in P^*$, and let $v\in S^*$ a different element.

    \noindent By Proposition \ref{prop:Omegaboxanddiamond} there are two sequences $(s_n)$, $(t_n)$ of elements of $S$ such that $\Omega_\Box(u,v)=\lim\limits_{m\rightarrow\infty}\lim\limits_{n\rightarrow\infty} \Omega(s_m,t_n)$ and \mbox{$\Omega_\Diamond(u,v)=\lim\limits_{n\rightarrow\infty}\lim\limits_{m\rightarrow\infty} \Omega(s_m,t_n).$}
    As $(s_n)$ is unbounded, 
    we have that $\Omega_\Diamond (u,v)=0.$

    \noindent Now, if $(t_n)$ is also unbounded, then $\lim\limits_{n\rightarrow \infty}\omega(t_n)=\infty$ too. Thus, $\Omega_\Box(u,v)=0$. If $(t_n)$ is bounded, then, for every $n,m \in \N$ (except maybe a finite number), $s_m\geq t_n$, and so 
        $$
            \Omega_\Diamond(u,v)=\lim_{m\rightarrow \infty} \lim_{n\rightarrow \infty} \frac{\omega(s_m \wedge t_n)}{\omega(s_n)\omega(t_n)}=\lim_{m\rightarrow \infty} \lim_{n\rightarrow \infty} \frac{1}{\omega(s_n)}=0.
        $$
    Thus
        \begin{align}\label{eq:Prod}
            \Omega_\Box(u,v)= \Omega_\Diamond (u,v)=0.
        \end{align}
    Thus, by (\ref{eq:Boxprod}), (\ref{eq:Diamondprod}) and (\ref{eq:Prod}), $\delta_u\in \mathfrak{Z}(\Aw'')$, but $\delta_u\notin \Aw$. Thus $\Aw$ is not strongly Arens irregular.

    An specific example of a weight that satisfies these conditions is $\omega(p)=1$, for $p\in [0,1]$ and $\omega(p)=p$, for $p>1.$
\end{ex}

\subsection{Approximate identities}
The existence of approximate identities is an interesting characteristic of Banach algebras. We shall see in the following subsection that the existence of approximate identities is intrinsically linked to duality and so to Arens regularity.

In the following result we refer to a net $(s_\alpha)$ of elements of $S$ that tends to $\sup S$. This means that, for any element $r\in S$, there exists $\alpha'$ such that $s_\alpha\geq r$ for any $\alpha\geq \alpha'.$

\begin{prop}\label{prop: Aw semigroup algebra BAI}
    Let $(S, \wedge)$ be an infinite, totally ordered semilattice and  $\omega$ a weight on $S$. Then, the weighted semigroup algebra $\Aw$ has an approximate identity. The following are true:
    \begin{itemize}
        \item[(a)] $\Aw$ has a bounded approximate identity if and only if there exists a net $(s_\nu)$ tending to $\sup S$ such that the set $\{\omega(s_\nu): \nu \}$ is bounded.
        \item[(b)] Suppose that for every strictly increasing net $(t_\nu)$ tending to $\sup S$, the set $\{\omega(t_\nu): \nu \}$ is unbounded. Then, $\Aw$ has a sequential approximate identity.
    \end{itemize} 
    $\Aw$ always has a multiplier-bounded approximate identity.
\end{prop}
\begin{proof}
    Let $\omega:S\rightarrow [1,\infty)$, and let $\alpha \in \Aw$. Then, for $s\in S$, we have
        \begin{align}\label{eq:alpha star delta sn}
            \alpha\star \delta_{s}=\sum_{t<s} \alpha(t)\delta_t + \sum_{t\geq s} \alpha(t)\delta_{s}
        \end{align}
        and
         \begin{align}\label{eq:alpha minus alpha star delta sn}
            \alpha-\alpha\star \delta_{s}= \sum_{t> s} \alpha(t)(\delta_t-\delta_{s}).
        \end{align}
    
    Suppose that there exists a strictly increasing net $(s_\nu)$ tending to $\sup S$ such that the set $\{\omega(s_\nu): \nu \}$ is bounded by $M\geq1$. Then, for $\alpha \in \Aw$, using (\ref{eq:alpha minus alpha star delta sn}) we see that
        \begin{align*}
            \|\alpha-\alpha\star \delta_{s_\nu}\|_\omega\leq \sum_{t> s_\nu} |\alpha(t)|\omega(t)+\left| \sum_{t > s_\nu} \alpha(t)\right| \omega(s_\nu)\leq (M+1) \sum_{t> s_\nu}  |\alpha(t)|\omega(t),
        \end{align*}
    which tends to zero since $\|\alpha\|_\omega<\infty$ and $s_\nu \rightarrow \sup S$.  Thus, $(s_\nu)$ is an approximate identity. Since $\| \delta_{s_\nu}\|_\omega=\omega(s_\nu)$, then it is a bounded approximate identity, with bound $M$. From (\ref{eq:alpha star delta sn}) we deduce that $(s_\nu)$ is a multiplier-bounded approximate identity. 
    
    On the other hand, suppose that for every strictly increasing net $(t_\nu)$ tending to $\sup S$, the set $\{\omega(t_\nu): \nu \}$ is unbounded. In this case, we can choose a strictly increasing sequence $(s_\nu)$ tending to $\sup S$ and such that $\omega(s_\nu)=\inf\{\omega(t): s_\nu\leq t\}$. To see this, consider the map $\tilde{\omega}: S \longrightarrow [1, \infty)$ such that \mbox{$\tilde{\omega} (s)=\inf\{ \omega(t): s\leq t \}$.} Since $\lim_{\nu} \omega(t_\nu)=\infty$ for every $(t_\nu)$ tending to $\sup S$, this infimum exists. Let $s_1$ such that $\omega(s_1)=\tilde{\omega}(1)$. Knowing $s_n$, take $s_{n+1}$ with $s_n<s_{n+1}$ and such that $ \omega(s_{n+1})=\tilde{\omega}(s_n+1)$. This sequence is such that $\omega(s_n)\leq w(t)$ for all $s_n\leq t$, and so, using (\ref{eq:alpha minus alpha star delta sn}), we have that
        \begin{align*}
            \|\alpha-\alpha\star \delta_{s_n}\|_\omega\leq \sum_{t> s_n} |\alpha(t)|\omega(t)+\left| \sum_{t > s_n} \alpha(t)\right| \omega(s_n)\leq 2 \sum_{t> s_n}  |\alpha(t)|\omega(t).
        \end{align*}
    Using the same reasoning as before, we see that $(s_n)$ defined this way is an approximate identity, and it is again a multiplier-bounded approximate identity. 
    
    \noindent Assume now towards contradiction that there is a bounded approximate identity in $\Aw$. Then $\|\cdot\|_\omega$ and $\|\cdot\|_{op}$ are equivalent. However 
        $$
            \lim_{n\rightarrow\infty}\|\delta_{s_n}\|_\omega=\lim_{n\rightarrow\infty}\omega(s_n)=\infty
        $$
     but $\|\delta_{s_n}\|_{op}\leq 1$, for all $n\in \N$. Thus, in this case $\Aw$ does not have a  bounded approximate identity, as desired.
\end{proof}

   Notice that, when $S$ is an infinite totally ordered semilattice such that $r=\sup S \in S$, then $\delta_r$ is an identity in $\Aw$, with $\| \delta_{r}\|_\omega=\omega(r)$. In fact we have a better result:
\begin{cor}\label{cor: identity for Aw}
  Let $S$ be a totally ordered semilattice. Then $\Aw$ has an identity if and only if $\sup S \in S$.
\end{cor}
\begin{proof}
    As seen above, if $\sup S \in S$, then $\Aw$ has an identity. Now suppose that $r=\sup S \notin S$ and assume towards contradiction that $\Aw$ has an identity $e$. Let us consider $(s_\nu)$ a net of elements in $S$ such that $(\delta_{s_\nu})$ is an approximate identity defined as in Proposition \ref{prop: Aw semigroup algebra BAI}. Recall that we can choose $(s_\nu)$ tending to $\sup S$. Then we have that
        \begin{align*}
            e=\lim_{\nu} \delta_{s_\nu}\star e= \lim_{\nu} \delta_{s_\nu}=\delta_r
        \end{align*}
    but we were supposing that $r\notin S$, and so $e=\delta_r \notin \Aw$.
\end{proof}

\begin{cor}\label{cor:S wedge Aw Tauberian}
    Let $(S, \wedge)$ be an infinite, totally ordered semilattice. Let $\omega:S\rightarrow [1,\infty)$ be a weight on $S$. Then when we view $\Aw$ as a Banach function algebra on $\varPhi_S$, it is Tauberian. 
\end{cor}
\begin{proof}
    This follows from the fact that the sequential approximate identity of $\Aw$ belongs to $c_{00}(S, 1/\omega)$.
\end{proof}

\subsection{(Non-) existence of Banach-algebra preduals}
As we have seen in Remark \ref{rem: linf dual Ew predual}, the space $E_\omega$ is a concrete predual of $\Aw$. However, it is not true that it is always a Banach-algebra predual. We shall study below in which situations we can know that $E_\omega$ is a predual of $\Aw$ and when it is unique. Also we shall give characterizations of when $\Aw$ is not a Banach-algebra predual for any concrete predual.

We shall start with some results that refer to generic semigroups, and we shall focus on totally ordered semilattices afterwards.

The following result is an extension of  \cite[Theorem 4.6]{DalesLauStrauss}, where they work only with the unweighted case.

\begin{prop}\label{theo:omegaboundeddual}
    Let $S$ be an infinite semigroup, and let $\omega$ be a weight on $S$. Suppose that $S$ is weakly cancellative. Then $E_\omega=c_0(S,1/\omega)$ is a Banach-algebra predual of $\Aw=\l1(S,\omega)$.

    Suppose that $S$ is not weakly cancellative. Suppose that there is a subset $U\subset S$ such that 
    \begin{enumerate}
        \item[(i)] $U$ is infinite;
        \item[(ii)]$\omega|U$ is bounded;
        \item[(iii)] there exist $t,u\in S$ such that $\{ r\in U: rt=u\} $ is infinite. 
    \end{enumerate}
    Then $E_\omega$ is not the predual of $\Aw$ as a Banach algebra.
\end{prop}
\begin{proof}
    Let $S$ be a weakly cancellative semigroup. Let $E_\omega=c_0(S, 1/\omega)$ and let $\lambda=(\lambda(s)) \in E_\omega$ and $\tilde{\delta_s}\in \Aw$. Then $\tilde{\delta_s}\cdot \lambda \in E_\omega$. Indeed, for $\varepsilon >0$ we have that there is a finite subset $F$ of $S$ such that 
        $$
            \left|\frac{\lambda(s)}{\omega(s)}\right|<\varepsilon \quad (s \in S \setminus F).
        $$
    Hence, as $S$ is weakly cancellative, the set $Fs^{-1}=\{r \in S: rs \in F \}$ is finite too. Thus
        \begin{align*}
            \left| \frac{(\tilde{\delta}_s\cdot \lambda)(r)}{\omega(r)}\right|&=\left|\frac{\lambda(rs)}{\omega(r)\omega(s)}\right|\leq\left|\frac{\lambda(rs)}{\omega(rs)}\right|<\varepsilon \quad (r\in S\setminus Fs^{-1}).
        \end{align*}
    And so $E_\omega$ is a submodule of $\Aw'$. Thus $\Aw$ is a dual Banach algebra as desired.

    Suppose now that $S$ is not weakly cancellative and that $U\subset S$ satisfies $(i)$, $(ii)$ and $(iii)$. Take $s_1,s_2,...\in U$ with $s_nt=u$ for every $n\in \N$.

    Consider $\alpha^{(n)}=\frac{1}{n}\left( \tilde{\delta}_{s_1}+\cdots+\tilde{\delta}_{s_n}\right)$ for $n\in\N$. Then $(\alpha^{(n)})$ is a sequence in $\Aw$ that tends to zero in $\sigma(\Aw,E_\omega)$. 

    Indeed, for $\lambda \in E_\omega$ and for $\varepsilon >0$ there exists a finite subset $F$ of $S$ such that $\left| \frac{\lambda(s)}{\omega(s)}\right| < \varepsilon$ for $s\in S \setminus F$. Thus $V=\{s_1,s_2,...\}\cap F$ is finite and so $\sum\limits_{s\in V}\left|\frac{\lambda(s)}{\omega(s)} \right|$ is bounded by $K>0$.
    Then
        \begin{align*}
            \left|\langle \alpha^{(n)}, \lambda \rangle \right| \leq \frac{1}{n}\sum\limits_{i=1}^{n} \left| \frac{\lambda(s_i)}{\omega(s_i)}\right|\leq \frac{K}{n}+\varepsilon.
        \end{align*}
    Thus $(\alpha^{(n)})$ tends to zero in $\sigma(\Aw,E_\omega).$

    Let $M<\infty$ be a bound for $\omega(s)$ ($s\in U$), and take $\lambda \in E_\omega$ such that $\langle\delta_u,\lambda\rangle\neq 0$. Then
        \begin{align*}
            \left|\langle\alpha^{(n)}\star \delta_t,\lambda\rangle\right|&
            \geq\frac{ 1}{M}\left|\langle\delta_u, \lambda\rangle\right|>0.
        \end{align*}
    Thus the multiplication in $\Aw$ is not separately $\sigma(\Aw,E_\omega)$-continuous, and so $E_\omega$ is not a submodule of $\Aw'$.
\end{proof}

As we have pointed out before, even some of the simplest examples of totally ordered semilattices are not even weakly cancellative, so they are not covered by Proposition \ref{theo:omegaboundeddual}. In the rest of this section we shall study what can be said for totally ordered semilattices.

\begin{prop}\label{theo:dualSomegaunbounded}
    Let $(S, \wedge)$ be an infinite, totally ordered semilattice. Let $\omega: S \longrightarrow [1,\infty)$ be a weight on $S$. Suppose that $\Lim\limits_{s\rightarrow \infty} \omega(s)=\infty$. Then $\Aw$ is a dual Banach algebra with Banach-algebra predual $E_\omega$. 
\end{prop}
\begin{proof}
    Let $\lambda=(\lambda(s)) \in E_\omega$ such that $\|\lambda\|_\omega'=1$, and let $\tilde{\delta}_t \in \Aw$ be the normalised point mass at $t\in S$ defined as above. Then, for $s\geq t$, we have
        \begin{align*}
            \frac{|(\tilde{\delta}_t\cdot \lambda)(s)|}{\omega(s)}=\frac{|\lambda(t)|}{\omega(s)\omega(t)} \leq \frac{1}{\omega(s)}. 
        \end{align*} 
    As in Definition \ref{def:limites}, for every $\varepsilon>0$, there is a finite subset $U_1$ of $S$ such that $\omega(s)>1/\varepsilon$ $(s\in S\setminus U_1)$. Hence in particular we have that $|(\tilde{\delta_t}\cdot \lambda)(s)|/\omega(s) <\varepsilon$ $(s\in  S\setminus U_1)$. 

    Then, for $s< t$, we have
        \begin{align*}
            \frac{|(\tilde{\delta}_t\cdot \lambda)(s)|}{\omega(s)}=\frac{|\lambda(s)|}{\omega(s)\omega(t)} \leq \frac{|\lambda(s)|}{\omega(s)}.
        \end{align*}  
    Since $\lambda \in E_\omega$, there exists $U_2$ a finite subset of $S$ such that $|\lambda(s)|/\omega(s) <\varepsilon$ for every $s\in S \setminus U_2$. 

    Thus $|(\tilde{\delta}_t\cdot \lambda)(s)|/\omega(s)<\varepsilon$ for every $S\setminus (U_1 \cup U_2)$. It follows that $(\tilde{\delta}_s\cdot \lambda) \in E_\omega$, and so $E_\omega$ is a closed submodule of $\Aw'$. Hence $\Aw$ is a dual Banach algebra as desired. 
\end{proof}

For an infinite semilattice (not necessarily totally ordered), the following result gives us a condition so that $E_\omega$ is unique as a Banach-algebra predual. Here we make use of results from Section 2 on character space of $\Aw$. We recall that, for a Banach algebra $A$, $L(A)=\mathrm{lin}\text{ }\varPhi_A$ and we see it as a linear subspace of $A'$.

\begin{prop}\label{prop:S abelian idempotent. E_w=L(character space) implies unique predual}
    Let $S$ be an infinite semilattice and let $\omega$ be a weight on $S$. Suppose that for all $x\in S$ the set $xS=\{xs: s\in S\}$ is finite and that $\Lim_{t\rightarrow \infty} \omega(t)=\infty$. Then:
    \begin{enumerate}
        \item[(i)] For every $x\in S$, $L_{\delta_x}: \Aw \rightarrow \Aw$ is compact.
        \item[(ii)] Suppose that $E_\omega=\overline{L(\Aw)}$. Then $E_\omega$ is the unique Banach-algebra predual of $\Aw$.
    \end{enumerate}
\end{prop}
\begin{proof}
    $(i)$ Let $x\in S$ and let $\varepsilon>0$. Let $M=\max \{ \omega(xt): t\in S\}$, which is well defined since $xS$ is finite. Since $\Lim_{t\rightarrow \infty} \omega(t)=\infty$, there exists $G$ a finite subset of $S$ such that $M/\omega(t)<\varepsilon$, $(t \in S\setminus G).$ Let $F=xG$, then 
        $$
            \{t\in S: xt \notin F\} \subseteq \{t \in S: t\notin G\}.
        $$

    \noindent Let $f \in \Aw$. Then, 
    \begin{align*}
         \|\pi_F \left( L_{\delta_x}(f)\right)- L_{\delta_x}(f)\|_\omega= \sum_{u\in S\setminus F} \left| \sum_{t\in S: xt=u} f(t)\right|\omega(u)\\
        \leq \sum_{u\in S\setminus F}  \sum_{t\in S: xt=u} |f(t)|\omega(u)= \sum_{t\in S:  xt\notin F} |f(t)|\omega(xt)\\
        \leq \|f\|_\omega \sup\{ \frac{M}{\omega(t)}: xt\notin F\}\leq \varepsilon \| f \|_\omega.
    \end{align*}

    \noindent Thus, since $\pi_F\circ L_{\delta_x}$ is finite rank and therefore compact, we obtain that $L_{\delta_x}$ is also compact.

    $(ii)$ Let $W$ be a Banach-algebra predual of $\Aw$. Assume towards contradiction that there exists $\varphi \in \varPhi_\omega \setminus W$. Then there exists $M\in \Aw''$ with $\|M\|_\omega\leq 1$ and such that $\langle M, \varphi\rangle=1$ and $\langle M, \lambda\rangle=0$ $(\lambda \in W)$. Thus there is a net $(f_\alpha)$ in $\Aw$ with norm bounded by 1 such that $\lim_\alpha f_\alpha(\varphi)=1$ and $\lim_\alpha \langle f_\alpha, \lambda\rangle=0$ $(\lambda \in W)$. Since $W$ is a Banach-algebra predual of $\Aw$, then we may suppose that there exists $f\in \Aw$ such that $\lim_\alpha f_\alpha=f$ in $\sigma(\Aw,W)$. Thus $\langle f,\lambda \rangle=0$ $(\lambda \in W)$. Hence $f=0$. Let $s\in S$ such that $\delta_s(\varphi)=1$, which exists by Proposition \ref{prop:character space of weighted semigroup algebras (idempotents)}. Since $\delta_s$ is compact we know that $\lim_\alpha \delta_s \star f_\alpha =\delta_s \star f$ in $\sigma(W',W'')=\sigma(A,A')$. Hence, $1=\lim_\alpha f_\alpha (\varphi)=f(\varphi)$, which is a contradiction, since $f\equiv 0$. Thus, $\varPhi_\omega \subset W$. But this implies that $\overline{L(\Aw)}\subset W$. Hence, by Proposition \ref{prop: concrete preduals subset}, $W=\overline{L(\Aw)}=E_\omega$ and so $E_\omega$ is the unique Banach-algebra predual of $\Aw$.
\end{proof}

In Proposition \ref{prop: N unique predual} we shall consider $\N$ with the minimum operation. Since it satisfies the conditions of Proposition \ref{prop:S abelian idempotent. E_w=L(character space) implies unique predual} that refer to the semigroup properties, we shall be able to use this proposition to identify for which weights the Banach-algebra predual $E_\omega$ is unique.

Under certain conditions, $\Aw$ is not a dual Banach algebra for any predual.

\begin{prop}\label{prop:dualomegaboundedSwedge}
    Let $(S, \wedge)$ be an infinite, totally ordered semilattice. Suppose that $\Liminf\limits_{s\rightarrow \infty} \omega(s) <\infty$. Suppose there exists an embedding of $S$ in $T$ specified as above, and that the set $\{s \in S: \omega(s)\leq \Liminf\omega \}$ has an accumulation point $r\in T\setminus S$. Then the Banach algebra $\Aw$ is not a dual Banach algebra with respect to any predual.
\end{prop}
\begin{proof}
    Let $U=\{s \in S: \omega(s)\leq \Liminf \omega \}$ and let us assume towards contradiction that there exists a Banach algebra predual $W$ for $\Aw$.

    Suppose that $r=\sup S$. Then, by Proposition \ref{prop: Aw semigroup algebra BAI}, $\Aw$ has a bounded approximate identity $(\delta_{s_\nu})$, with $(s_\nu)$ tending to $r$. Since $\Aw$ is a dual Banach algebra, there exists a subnet $(\delta_\alpha)$ converging in the topology $\sigma(\Aw,W)$ to an identity $e\in \Aw$, but this is a contradiction with Corollary \ref{cor: identity for Aw}.

    Suppose that $r\neq \sup S$. Then there exists a net $(s_\beta)$ in $U$ monotone decreasing or increasing converging to $r$. Suppose that $(s_\beta)$ is decreasing. Since $\|\delta_{s_\beta}\|_\omega$ is bounded for every $\beta$, there exists a subnet $(\delta_\alpha)$ that converges to an element $f\in \Aw$ in the topology $\sigma(A,W)$. We know that $\|f\|_\omega=\lim_{\alpha}\|\delta_\alpha\|_\omega\geq 1$. Thus, $f$ is a non-zero element of $\Aw$. Let $s\in S$ such that $s>r$. Then $ \delta_s \star \delta_\alpha=\delta_\alpha$ for every $\alpha$ large enough. Since the multiplication is separately $\sigma(\Aw,W)$-continuous $\delta_s\star f=f$. This implies that $\supp f \subset (0,s]\cap S$ for every $s>r$. Thus $\supp f \subset (0,r]\cap S$. Now let $s\in S$ such that $s<r$. Then $\delta_s\star \delta_\alpha=\delta_s$. Thus $\delta_s\star f=\delta_s$, which implies that $f$ is not zero and that $\supp f \subset [s,\sup S)\cap S$. Since this is true for every $s<r$ then $\supp f \subset [r,\sup S)\cap S$. We conclude that $\supp f \subset \{r\}\cap S$. But $\{r\}\cap S =\emptyset$, and so there is no such $f\in \Aw$. The case where $(s_\beta)$ is increasing is symmetrical. Thus $\Aw$ is not a dual Banach algebra. 
\end{proof}

\begin{ex}\label{ex:Qnotinf}
    Let us consider again the semigroup and weight defined in Example \ref{ex:QnotSAI}. 
    Since any positive irrational number is an accumulation point of this set, by Proposition \ref{prop:dualomegaboundedSwedge} $\Aw$ is not a dual Banach algebra.
\end{ex}

\section{Specific results for $(\N,\wedge)$}\label{sec:N_wedge}

\subsection{Arens Regularity and strong Arens irregularity, DTC sets}

Consider the semigroup $S:= \N$ with the semigroup operation
    $$
        m\wedge n=\min\{ m,n\} \quad (m,n \in \N),
    $$
which is a particular case of the above. Throughout this section we shall write $D_\omega=\ell^1(\N_\wedge, \omega)$, where $\omega$ is a weight throughout this section.

We shall give below some results that improve what we have obtained in the more generic case.

In the previous section we have seen when a weighted semigroup algebra is not strongly Arens irregular, however we do not have a characterization of when it is strongly Arens irregular. In the following result we shall see that for $S=(\N, \wedge)$, not only we can see when $D_\omega$ is strongly Arens irregular, but we can also determine the smallest DTC set. In \cite[Example 7.33]{DalesLauStrauss} they work with the unweighted case and they prove that it is strongly Arens irregular, which also follows from \cite[Theorem 2.14]{DalesStrauss}. The argument followed here is very similar to the one followed in \cite{DalesLauStrauss}.

\begin{deff}
    Let $A$ be a Banach algebra. Then a subset $V$ of $A''$ is \textit{ determining for the left topological centre} (a DLTC set) of $A''$ if, given $M\in A''$ such that $M\Box N= M\Diamond N$ $(N\in V)$, then $M\in A$. When the algebra $A$ is commutative, we use \textit{ determining for the topological centre} (a DTC set).
\end{deff}

\begin{prop}\label{cor:DTCsetNminsemigroupalg}
    Let $S=\N$ with the semigroup operation $\wedge$ defined as above. Then whenever $\liminf_{n\rightarrow \infty}\omega(n) <\infty$ the Banach algebra $D_\omega$ is strongly Arens irregular. Furthermore, there is a two-point DTC set of $D_\omega''$.
\end{prop}
\begin{proof}Consider the isometric isomorphism $\theta_\omega$ defined in Lemma \ref{lem: theta_omega}.

    We have that
        \begin{align*}
            \delta_u\Box \delta_v =\delta_u,& \quad (u\in \beta \N, v\in \N^*);\\
            \delta_u\Diamond \delta_v=\delta_v,& \quad (u\in \N^*, v\in \beta\N).
        \end{align*}

    Let $u,v \in \beta \N$ and $(s_\alpha)$, $(t_\beta)$ nets in $\N$ such that $u=\lim_\alpha s_\alpha$ and $v=\lim_\beta t_\beta$. Thus, for $\lambda \in C(\beta \N)$, we have
        \begin{align*}
            \langle \theta''_\omega(\delta_u)\Box_\omega \theta''_\omega(\delta_v), \omega \lambda \rangle =& \lim_\alpha \lim_\beta  \theta''_\omega(\delta_{s_\alpha})\Box_\omega \theta''_\omega(\delta_{t_\beta}), \omega \lambda \rangle\\
            &\Omega_\Box (u,v) \langle \delta_u \Box \delta_v, \lambda\rangle = \Omega_\Box(u,v)\langle \delta_u,\lambda\rangle.
        \end{align*}
    Symmetrically we obtain that
        \begin{align*}
            \langle \theta''_\omega(\delta_u)\Diamond_\omega \theta''_\omega(\delta_v), \omega \lambda \rangle = \Omega_\Diamond(u,v)\langle \delta_u,\lambda\rangle.
        \end{align*}

    Consider $\varphi_\omega \in D_\omega'$ defined as $\varphi_\omega(\alpha)=\sum\limits_{n=1}^\infty \alpha(n)/\omega(n)$, $(\alpha\in D_\omega)$. Let $u\in\N^*$ and say $u=\lim_\alpha s_\alpha$, where $(s_\alpha)$ is a net in $\N$. Then, for $t\in \N$, we have
        \begin{align*}
            \Omega_\Box(t,u)=\lim_\alpha 1/\omega(s_\alpha)=\lim_\alpha\langle \delta_{s_\alpha}, \varphi_\omega\rangle = \langle \delta_u,\varphi_\omega\rangle.
        \end{align*}
    Let $u\in\N$, $g \in M(\N^*)$ and $\lambda \in C(\beta \N).$ Then
         \begin{align*}
            \langle \theta''_\omega(\delta_u)\Box_\omega \theta''_\omega(g), \omega \lambda\rangle =& \langle \sum\limits_{s\in \N^*} g(s)\theta''_\omega(\delta_u)\Box_\omega \theta''_\omega(\delta_s),\omega \lambda\rangle\\
                =\langle \Omega_\Box(u,s)g(s)\delta_u,\lambda\rangle\\
                =\langle g,\varphi_\omega\rangle \langle\delta_u,\lambda\rangle.
        \end{align*}
    Thus, for $\mu \in M(\beta \N), \nu \in M(\N^*)$, we have
        \begin{align*}
            \theta''_\omega(\mu)\Box_\omega\theta''_\omega(\nu)=\langle \nu, \varphi_\omega\rangle \mu.
        \end{align*}
    Symmetrically, we obtain
        \begin{align*}
            \theta''_\omega(\nu)\Diamond_\omega\theta''_\omega(\mu)=\langle \mu, \varphi_\omega\rangle\nu \quad(\nu \in M(\beta \N), \mu \in M(\N^*))
        \end{align*}
    Hence $\theta''_\omega(M(\N^*))$ is a closed subalgebra of $D_\omega''$. Also for $\mu \in M(\beta\N)$ and $f\in \ell^1(\N)$, we have
        \begin{align*}
            \theta''_\omega (f) \cdot \theta''_\omega(\mu)=\theta''_\omega(\mu)\cdot \theta''_\omega(f)=\langle \mu,\varphi_\omega\rangle f.
        \end{align*}
    Hence $D_\omega$ is an ideal in its bidual and we can write 
        $$
            D_\omega''=\theta''_\omega(M(\N^*))\ltimes D_\omega.
        $$

    Let $a,b \in \N^*$ be two different points such that  at least one of the quantities $\langle \delta_a,\varphi_\omega\rangle$ or $\langle \delta_b,\varphi_\omega \rangle$ is not zero, which exist since $\Lim\limits_{s\rightarrow \infty}\omega(s)<\infty$. Take \mbox{$\mu \in M(\N^*)$} such that  
        $$
        \theta''_\omega (\mu)\Box_\omega \theta''_\omega (\delta_a)= \theta''_\omega (\mu)\Diamond_\omega \theta''_\omega (\delta_a) \text{ and } \theta''_\omega (\mu)\Box_\omega \theta''_\omega (\delta_b)= \theta''_\omega (\mu)\Diamond_\omega \theta''_\omega (\delta_b).
        $$ 
    Then $\theta''_\omega(\mu)=0$. Thus $D_\omega$ is strongly Arens irregular and \mbox{$V=\{\theta''_\omega(\delta_a),\theta''_\omega(\delta_b)\}$} is a DTC set for $D_\omega''$.  
\end{proof}

\begin{rem}
    In \cite[Example 9.13]{Dedania}  it is observed that $D_\omega$ is Arens regular for any weight $\omega$ such that $\omega(n)\rightarrow\infty$ when $n\rightarrow \infty$. However, the authors justify this by appealing to a theorem only stated for cancellative semigroups, while $\N_\wedge$ is not even weakly cancellative. By applying our Theorem \ref{theo:Arens regularity when Limit of omega is infinity} the desired result follows easily.
\end{rem}

To finish, we shall look at the duality of $D_\omega$.

    \begin{prop}\label{prop: N unique predual}
    Let $\omega:\N\rightarrow [1,\infty)$. Then $D_\omega$ is a dual Banach algebra if and only if $\lim_{n\rightarrow \infty} \omega(n)=\infty$. In this case, $E_\omega$ is the unique Banach-algebra predual.
    \end{prop}
    \begin{proof}
        When  $\liminf\limits_{n\in \infty} \omega(n) <\infty$ the Banach algebra $D_\omega$ is not a dual Banach algebra, as follows from Proposition \ref{prop:dualomegaboundedSwedge}. 
    
        It follows from Proposition \ref{theo:dualSomegaunbounded} that $D_\omega$ is a dual Banach algebra with predual $E_\omega$ when $\omega: \N\longrightarrow [1,\infty)$ is a weight on $\N_\wedge$ such that $\lim\limits_{n\rightarrow \infty} \omega(n)=\infty$. What is more, let $\varphi\in \varPhi_\omega$. Then there exists $k\in \N$ such that $\varphi=\varphi_k$ where
            \begin{equation}\label{eq:varphi_k}
                \varphi_k(\alpha)=\sum_{n=k}^\infty \alpha(n) \quad (\alpha=(\alpha(n))\in D_\omega).
            \end{equation}
        When $\lim_{n\rightarrow \infty} \omega(n)=\infty$, we have that $\varphi_\omega \in E_\omega$ and so $\varPhi_\omega \subset E_\omega$. Hence $\overline{L(D_\omega)}\subset E_\omega$.
    
        \noindent Take $\lambda=(\lambda(n))\in E_\omega$. For $k\in \N$, consider $\rho^{(k)}$ defined as follows:
        \begin{align*}
            \rho^{(k)}= \varphi_k \lambda(k)-\varphi_{k+1}\lambda(k),
        \end{align*}
        where $\varphi_k$ is defined in (\ref{eq:varphi_k}). Then, since $\lambda \in E_\omega$, we have that $\rho^{(k)} \in \varPhi_\omega$ and so $\mu^{(n)} =\sum_{i=1}^k \rho^{(i)} \in L(D_\omega)$, for $n\in \N$. Finally, 
        $$
            \|\lambda-\mu^{(n)}\|_\omega= \|\sum_{i=n}^{\infty} \lambda(i)\|_\omega,
         $$
         since $\lambda \in E_\omega$, $\|\lambda-\mu^{(n)}\|_\omega \rightarrow 0$ when $n\rightarrow \infty.$ Thus, $\lambda \in \overline{L(D_\omega)}$. So, we can apply Proposition \ref{prop:S abelian idempotent. E_w=L(character space) implies unique predual} to obtain that $E_\omega$ is the unique Banach-algebra predual of $D_\omega.$
    \end{proof}

To conclude the study of $D_\omega,$ notice that by Proposition \ref{prop: Aw semigroup algebra BAI}, we obtain that $D_\omega$ has a bounded approximate identity if and only if \mbox{$\liminf\limits_{n\rightarrow \infty} \omega(n) < \infty$} and that it always has a multiplier-bounded approximate identity.

\section*{Acknowledgement}
To my supervisor, Professor H. G. Dales, in memoriam, for having shared with me his passion for Mathematics until the last moment.

\renewcommand{\bibname}{References}

\printbibliography

@article {Arens-2,
    AUTHOR = {Arens, Richard},
     TITLE = {The adjoint of a bilinear operation},
   JOURNAL = {Proc. Amer. Math. Soc.},
  FJOURNAL = {Proceedings of the American Mathematical Society},
    VOLUME = {2},
      YEAR = {1951},
     PAGES = {839--848},
      ISSN = {0002-9939},
   MRCLASS = {46.3X},
  MRNUMBER = {45941},
MRREVIEWER = {L. W. Cohen},
       DOI = {10.2307/2031695},
       URL = { },
}

@article {Arens-1,
    AUTHOR = {Arens, Richard},
     TITLE = {Operations induced in function classes},
   JOURNAL = {Monatsh. Math.},
  FJOURNAL = {Monatshefte f\"{u}r Mathematik},
    VOLUME = {55},
      YEAR = {1951},
     PAGES = {1--19},
      ISSN = {0026-9255},
   MRCLASS = {56.0X},
  MRNUMBER = {44109},
MRREVIEWER = {L. W. Cohen},
       DOI = {10.1007/BF01300644},
       URL = { },
}

@article {choi_2013,
    AUTHOR = {Choi, Yemon},
     TITLE = {Approximately multiplicative maps from weighted semilattice
              algebras},
   JOURNAL = {J. Aust. Math. Soc.},
  FJOURNAL = {Journal of the Australian Mathematical Society},
    VOLUME = {95},
      YEAR = {2013},
    NUMBER = {1},
     PAGES = {36--67},
      ISSN = {1446-7887},
   MRCLASS = {46J10 (39B72)},
  MRNUMBER = {3123743},
MRREVIEWER = {Miguel Cabrera},
       DOI = {10.1017/S1446788713000189},
       URL = { },
}

@article {CrawYoung,
    AUTHOR = {Craw, I. G. and Young, N. J.},
     TITLE = {Regularity of multiplication in weighted group and semigroup
              algebras},
   JOURNAL = {Quart. J. Math. Oxford Ser. (2)},
  FJOURNAL = {The Quarterly Journal of Mathematics. Oxford. Second Series},
    VOLUME = {25},
      YEAR = {1974},
     PAGES = {351--358},
      ISSN = {0033-5606},
   MRCLASS = {43A20},
  MRNUMBER = {365029},
MRREVIEWER = {G. P. Johnson},
       DOI = {10.1093/qmath/25.1.351},
       URL = { },
}

@article {Dedania,
    AUTHOR = {Dales, Harold Garth and Dedania, Haresh V.},
     TITLE = {Weighted convolution algebras on subsemigroups of the real
              line},
   JOURNAL = {Dissertationes Math.},
  FJOURNAL = {Dissertationes Mathematicae},
    VOLUME = {459},
      YEAR = {2009},
     PAGES = {60},
      ISSN = {0012-3862},
   MRCLASS = {46H20 (43A20)},
  MRNUMBER = {2477218},
MRREVIEWER = {Ali \"{U}lger},
       DOI = {10.4064/dm459-0-1},
       URL = { },
}

@article {DalesLau,
    AUTHOR = {Dales, H. G. and Lau, A. T.-M.},
     TITLE = {The Second Duals of {B}eurling Algebras},
   JOURNAL = {Mem. Amer. Math. Soc.},
  FJOURNAL = {Memoirs of the American Mathematical Society},
    VOLUME = {177},
      YEAR = {2005},
    NUMBER = {836},
     PAGES = {vi+191},
      ISSN = {0065-9266},
   MRCLASS = {43A10 (43A20 46J10)},
  MRNUMBER = {2155972},
MRREVIEWER = {Ali \"{U}lger},
       DOI = {10.1090/memo/0836},
       URL = { },
}

@article {DalesLoy,
    AUTHOR = {Dales, H. G. and Loy, R. J.},
     TITLE = {Approximate Amenability of Semigroup Algebras and {S}egal
              Algebras},
   JOURNAL = {Dissertationes Math.},
  FJOURNAL = {Dissertationes Mathematicae},
    VOLUME = {474},
      YEAR = {2010},
     PAGES = {58},
      ISSN = {0012-3862},
   MRCLASS = {46H20 (43A20 46M20)},
  MRNUMBER = {2760649},
MRREVIEWER = {George A. Willis},
       DOI = {10.4064/dm474-0-1},
       URL = { },
}

@article {DalesLauStrauss,
    AUTHOR = {Dales, H. G. and Lau, A. T.-M. and Strauss, D.},
     TITLE = {Banach Algebras on Semigroups and on their Compactifications},
   JOURNAL = {Mem. Amer. Math. Soc.},
  FJOURNAL = {Memoirs of the American Mathematical Society},
    VOLUME = {205},
      YEAR = {2010},
    NUMBER = {966},
     PAGES = {vi+165},
      ISSN = {0065-9266},
      ISBN = {978-0-8218-4775-6},
   MRCLASS = {43A10 (43A20 46H05)},
  MRNUMBER = {2650729},
MRREVIEWER = {Volker Runde},
       DOI = {10.1090/S0065-9266-10-00595-8},
       URL = { },
}

@article {DalesStrauss,
    AUTHOR = {Dales, H. G. and Strauss, D.},
     TITLE = {Arens regularity for totally ordered semigroups},
   JOURNAL = {Semigroup Forum},
  FJOURNAL = {Semigroup Forum},
    VOLUME = {105},
      YEAR = {2022},
    NUMBER = {1},
     PAGES = {172--190},
      ISSN = {0037-1912},
   MRCLASS = {Prelim},
  MRNUMBER = {4466963},
       DOI = {10.1007/s00233-022-10299-x},
       URL = { },
}

@book{DalesUlger,
  title={Banach Function Algebras, Arens Regularity, and BSE Norms},
  author={H. G. Dales and A. Ulger},
  year={In preparation},
  publisher={}
}

@article {Daws-Connes-Amenability,
    AUTHOR = {Daws, Matthew},
     TITLE = {Connes-amenability of bidual and weighted semigroup algebras},
   JOURNAL = {Math. Scand.},
  FJOURNAL = {Mathematica Scandinavica},
    VOLUME = {99},
      YEAR = {2006},
    NUMBER = {2},
     PAGES = {217--246},
      ISSN = {0025-5521},
   MRCLASS = {46H05 (43A07 46B10)},
  MRNUMBER = {2289023},
MRREVIEWER = {Volker Runde},
       DOI = {10.7146/math.scand.a-15010},
       URL = { },
}

@article {Feinstein,
    AUTHOR = {Feinstein, J. F.},
     TITLE = {Strong {D}itkin algebras without bounded relative units},
   JOURNAL = {Int. J. Math. Math. Sci.},
  FJOURNAL = {International Journal of Mathematics and Mathematical
              Sciences},
    VOLUME = {22},
      YEAR = {1999},
    NUMBER = {2},
     PAGES = {437--443},
      ISSN = {0161-1712},
   MRCLASS = {46J10},
  MRNUMBER = {1695574},
MRREVIEWER = {Pamela Gorkin},
       DOI = {10.1155/S0161171299224374},
       URL = { },
}

@article {filalitopcent,
    AUTHOR = {Filali, Mahmoud and Salmi, Pekka},
     TITLE = {Topological centres of weighted convolution algebras},
   JOURNAL = {J. Funct. Anal.},
  FJOURNAL = {Journal of Functional Analysis},
    VOLUME = {278},
      YEAR = {2020},
    NUMBER = {11},
     PAGES = {108468, 22},
      ISSN = {0022-1236},
   MRCLASS = {43A20 (22D15)},
  MRNUMBER = {4075579},
MRREVIEWER = {Abdolrasoul Pourabbas},
       DOI = {10.1016/j.jfa.2020.108468},
       URL = { },
}

@article {Lust-piquard,
    AUTHOR = {Lust-Piquard, Fran\c{c}oise},
     TITLE = {\'{E}l\'{e}ments ergodiques et totalement ergodiques dans
              {$L^{\infty }(\Gamma )$}},
   JOURNAL = {Studia Math.},
  FJOURNAL = {Polska Akademia Nauk. Instytut Matematyczny. Studia
              Mathematica},
    VOLUME = {69},
      YEAR = {1980/81},
    NUMBER = {3},
     PAGES = {191--225},
      ISSN = {0039-3223},
   MRCLASS = {43A07 (43A60)},
  MRNUMBER = {647138},
MRREVIEWER = {S. Hartman},
       DOI = {10.4064/sm-69-3-191-225},
       URL = { },
}

@article {Ross,
    AUTHOR = {Ross, Kenneth A.},
    TITLE = {The structure of certain measure algebras},
    JOURNAL = {Pacific J. Math.},
    FJOURNAL = {Pacific Journal of Mathematics},
    VOLUME = {11},
    YEAR = {1961},
    PAGES = {723--737},
    ISSN = {0030-8730},
    MRCLASS = {46.80 (42.56)},
    MRNUMBER = {137003},
    MRREVIEWER = {D. A. Edwards},
       URL = { },
}

@article {White,
    AUTHOR = {White, Katherine},
     TITLE = {Amenability and ideal structure of some {B}anach sequence
              algebras},
   JOURNAL = {J. London Math. Soc. (2)},
  FJOURNAL = {Journal of the London Mathematical Society. Second Series},
    VOLUME = {68},
      YEAR = {2003},
    NUMBER = {2},
     PAGES = {444--460},
      ISSN = {0024-6107},
   MRCLASS = {46J20},
  MRNUMBER = {1994693},
MRREVIEWER = {Joel Francis Feinstein},
       DOI = {10.1112/S0024610703004575},
       URL = { },
}

\end{document}